\documentclass[12pt]{amsart}
\usepackage{amsmath,amstext,amsfonts,amsthm,amssymb}
\usepackage{eucal}
\usepackage{diagrams}
\diagramstyle[scriptlabels,PostScript=dvips]
\renewcommand{\subsubsection}[1]{\addtocounter{subsubsection}{1}
{\ \\[3pt]\bf \thesubsubsection. \  #1} }
\swapnumbers

\newtheorem{thm}[subsubsection]{Theorem}
\newtheorem{lem}[subsubsection]{Lemma}
\newtheorem{prp}[subsubsection]{Proposition}

{  \theoremstyle{definition}
           \newtheorem{dfn}[subsubsection]{Definition}
           \newtheorem{rems}[subsubsection]{Remarks}

}

\newcommand{\one}{\mathbf{1}}       
\newcommand{\Sh}{\mathrm{Sh}}       
\newcommand{\DD}{\cD}        

\newcommand{\et}{\mathrm{et}}
\newcommand{\ET}{\mathrm{ET}}
\newcommand{\Mod}{\mathrm{Mod}}

\renewcommand{\H}{\CMcal{H}}
\newcommand{\cC}{\mathcal{C}}
\newcommand{\cD}{\mathcal{D}}

\newcommand{\cM}{\mathcal{M}}
\newcommand{\cO}{\mathcal{O}}

\newcommand{\cU}{\mathcal{U}}

\newcommand{\cX}{\mathcal{X}}
\newcommand{\cY}{\mathcal{Y}}

\newcommand{\id}{\mathrm{id}}

\newcommand{\C}{\mathbb{C}}

\newcommand{\isom}{\approx}

\newcommand{\Hom}{\operatorname{Hom}}

\newcommand{\Ker}{\mathrm{Ker}}

\newcommand{\op}{\mathrm{op}}

\newcommand{\str}{\mathrm{str}}
\newcommand{\Aut}{\mathrm{Aut}}
\newcommand{\Atlas}{\mathtt{Atlas}}
\newcommand{\bs}{\backslash}
\newcommand{\coc}{\mathrm{coc}} 
\newcommand{\Ob}{\operatorname{Ob}}

\begin{document}

%
%

\title[]{Drinfeld double for orbifolds}
\author{Vladimir Hinich}    
\address{Department of Mathematics, University of Haifa,    
Mount Carmel, Haifa 31905,  Israel}    
\email{hinich@math.haifa.ac.il}    
\begin{abstract}
We prove that the Drinfeld double of the category of sheaves
on an orbifold is equivalent the the category of sheaves on the corresponding
inertia orbifold. We hope this observation will help to explain the
appearence of inertia groupoid in the orbifold cohomology theory.
\end{abstract}  
\maketitle    

{\center{\it Dedicated to the memory of Iosif Donin}}

\section{Introduction}
\subsection{Overview}

This note is a result of an attempt to understand why inertia orbifold 
appears in the stringy cohomology theory. Recall that if $\cX$ is an orbifold
then its (stringy) orbifold cohomology is defined by the formula
\begin{equation}
\H_\str(\cX)=\H(I\cX)
\end{equation} 
where $\H$ is a usual cohomology for orbifolds and $I\cX$ is the inertia
orbifold which is (roughly speaking) the collection of pairs $(x,\gamma)$
with $x\in\cX$ and $\gamma\in\Aut(x)$.

In an attempt to find a ``cohomological meaning'' of the connection 
between $\cX$ and $I\cX$ we have found a very simple way to express
the category of sheaves on $I\cX$ through that on $\cX$.

\begin{thm}The category of sheaves of $\cO$-modules $\Sh(I\cX)$ is
naturally equivalent to the Drinfeld double of the category $\Sh(\cX)$
(considered as a monoidal category with respect to the tensor product
of sheaves).
\end{thm}

The orbifolds are assumed to be complex or $C^\infty$-orbifolds. In the
complex case one can substitute the sheaves of $\cO$-modules with the
coherent sheaves.

The functor $\Theta:\Sh(I\cX)\rTo\DD(\Sh(\cX))$ from the category of sheaves on
$I\cX$ to the Drinfeld double of $\Sh(\cX)$ is constructed as follows.

The canonical projection $\pi:I\cX\rTo\cX$ is defined by the formula
$\pi(x,\gamma)=x$. This is a finite morphism of orbifolds and 
one has the projection formula
\begin{equation}\label{intro:proj}
\pi_*(M\otimes \pi^*(N))=\pi_*(M)\otimes N.
\end{equation}
For any $N\in\Sh(\cX)$ the inverse image $\pi^*(N)$ admits a canonical 
automorphism $\zeta_N:p^*(N)\rTo p^*(N)$ 
defined at each point $(x,\gamma)\in I\cX$ by the formula
$$ (\pi^*(N))_{(x,\gamma)}=N_x\rTo^\gamma N_x=(\pi^*(N))_{(x,\gamma)}.$$
Thus, for any sheaf $M\in\Sh(I\cX)$ the functor 
$N\mapsto \pi_*(M)\otimes N=\pi_*(M\otimes \pi^*(N))$
is endowed with a canonical endomorphism induced by $\zeta_N$. Its
composition with the standard commutativity constraint
gives an object of $\DD(\Sh(\cX))$.

To prove $\Theta$ is an equivalence, we check this locally (for a quotient
of a manifold by a finite group) and then use the local structure of
complex or $C^\infty$-orbifolds to deduce the general claim.
\footnote{
We believe our result is also valid for DM stacks (and for quasicohorent 
or coherent sheaves). However, our proof
which relies upon the existence of an open covering of the orbifold
with the global quotients, does not work in the algebraic case.
}

Note that the equivalence $\Theta$ provides the category $\Sh(I\cX)$
with a new structure of braded monoidal category. We explicitly describe
this new monoidal structure in~\ref{newbradedstructure}.

\subsection{Structure of the paper}
In Section~\ref{general} we recall general facts about orbifolds,
including the description of the standard functors.
In Section~\ref{theta} we recall the definition of Drinfeld double of a 
monoidal category and construct the functor $\Theta$.
In Sections~\ref{affine} and \ref{affine2} we prove that $\Theta$ is
an equivalence in the affine case: for a quotient of an affine scheme by a 
finite group or for a (global) quotient of a manifold (resp., Stein manifold)
by a finite group.
In Section~\ref{sec:atlas} we describe a presentation of an orbifold
by an atlas of affine orbifold charts.
This allows us to get the general result using an affine orbifold atlas
in Section~\ref{localtoglobal}.

\section{Generalities on orbifolds}\label{general}

We define an orbifold as a stack on a fixed category $\cM$ of manifolds
having a ``geometric origin''. It makes sense considering the following
instances of the base category $\cM$.
\begin{itemize}
\item $C^\infty$ manifolds.
\item Complex manifolds.
\item Schemes over a fixed base scheme $S$.\footnote{
and, probably, many others.}
\end{itemize}
The category $\cM$ of one of the above mentioned types is
endowed with the \'etale topology: the coverings are finite collections 
of \'etale (=local isomorphisms in the $C^\infty$ or complex case) morphisms.

\subsection{Groupoids}
An \'etale groupoid $X_\bullet$ in $\cM$ is a collection of the following data.
\begin{itemize}
\item A space $X_0\in\cM$ of objects.
\item A space $X_1\in\cM$ of arrows.
\item \'Etale morphisms of source and target $s,t:X_1\rTo X_0$.
\item Composition map $X_1\times_{X_0}X_1\rTo X_1$.
\item An inversion map $i:X_1\rTo X_1$. 
\end{itemize}
The requirements imposed on these data guarantee that for each $M\in\cM$
the data $\Hom(M,X_\bullet)$ form a groupoid.

Note that since the source map is \'etale, the fibre product
$X_1\times_{X_0}X_1$ exists in $\cM$. A groupoid $X$ is called separated
if the map
$$(s,t):X_1\rTo X_0\times X_0$$
is proper.

If, for example, $X$ is a manifold and $G$ is a finite group
acting on $X$, the quotient groupoid $(G\bs X)_\bullet$ is defined by
the formulas
$$ (G\bs X)_0=X,\quad  (G\bs X)_1=G\times X,\ s(g,x)=x,\ t(g,x)=g(x).$$

Each groupoid $X_\bullet$ in $\cM$ defines a presheaf of groupoids
on $\cM$ assigning to $M\in\cM$ a groupoid $\Hom(M,X_\bullet)$.

The associated stack of groupoids is denoted $[X_\bullet]$.
In particular, $[G\bs X]$ denotes the stack associated to $(G\bs X)_\bullet$.

\subsection{Orbifolds}

\begin{dfn}
An orbifold $\cX$ is a stack of groupoids on $\cM$
equivalent to a stack of the form $[X_\bullet]$ where
$X_\bullet$ is an \'etale separated groupoid.
\end{dfn}

\subsection{Affine orbifolds}
Let $G$ be a finite group acting on an affine smooth algebraic variety $X$.
The orbifold equivalent to $[G\bs X]$ is called {\em an affine orbifold}.

Similarly, in the complex-analytic  version, 
$X$ is assumed to be a Stein manifold; in the $C^\infty$ version
$X$ can be any smooth manifold.

An important property of complex-analytic and $C^\infty$-orbifolds is the
existence of affine open neighborhoods.

\begin{prp}
\label{prp:orbi-neighborhood}
Let $\cM$ be the category of complex or $C^\infty$ manifolds.
Let $X_\bullet$ be an \'etale separated groupoid in $\cM$.
Let $x\in X_0$ and $G=\Aut(x)=\{\gamma\in X_1|s(\gamma)=t(\gamma)=x\}$.
There exists an open neighborhood $U$ of $x$ in $X_0$ (Stein in the
complex case) so that
the restriction of $X_\bullet$ to $U$  is isomorphic to a quotient
groupoid $(G\bs U)_\bullet$.
\end{prp}
\begin{proof}
The $C^\infty$ case is proven in~\cite{moep}, proof of Th.~4.1.
In the complex case one should only add that $U$ is obtained as a finite
intersection of open polydiscs, and is, therefore, Stein.
\end{proof}

\subsection{Inertia orbifold}

Let $G$ be a groupoid. We define the groupoid $IG$ (inertia groupoid of $G$)
as follows. The objects of $IG$ are pairs $(x,\gamma)$ with $x\in G,
\gamma:x\rTo x$. A morphism from $(x,\gamma)$ to $(y,\delta)$ is a morphism
$\alpha:x\rTo y$ satisfying the property $\delta\alpha=\alpha\gamma$.

If $\cX$ is an orbifold, we define $I\cX$ as the fibered category on $\cM$
assigning a groupoid $I(\cX(M))$ to $M\in\cM$.

\begin{lem}
$I\cX$ is an orbifold.
\end{lem}
\begin{proof}
Let us check first that $I\cX$ is a stack. Let 
$(x,\gamma),(y,\delta)\in I\cX(M)$. One has
$$\Hom((x,\gamma),(y,\delta))=\Ker\left(\Hom(x,y)\pile{\rTo^{\gamma^*}\\
\rTo_{\delta_*}}\Hom(x,y)\right)$$
and this proves that Hom-sets in $I\cX$ form a sheaf.
The second axiom of stacks saying that the groupoid $I\cX(U),\ U\in\cM,$
is equivalent to the groupoid of the descent data with respect to any covering
of $U$, is straightforward.

Finally, suppose $\cX=[X_\bullet]$. We claim that $I\cX$ is equivalent to
$[IX_\bullet]$ where the groupoid $IX_\bullet$ is defined by the formulas
\begin{eqnarray}
IX_0&=&\Ker(X_1\pile{\rTo^s\\ \rTo_t}X_0)\\
IX_1&=&\Ker(X_1\times X_1\pile{\rTo^{s\times t}\\ \rTo_{t\times s}}
X_0\times X_0).
\end{eqnarray}
An important step in proving the claim is to check that $IX_0$ and $IX_1$
are smooth manifolds. 
In complex or $C^\infty$ case this follows from
Proposition~\ref{prp:orbi-neighborhood}. In algebraic case the groupoid
$IX_\bullet$ as above is defined; its smoothness can be checked locally
in \'etale topology, and this also reduces the claim to the affine case.
\end{proof}

\subsection{Sheaves on orbifolds}

Let $\cX$ be an orbifold.
We define the \'etale cite of $\cX$, $\cX_\et$, as follows.

The objects of $\cX_\et$ are \'etale morphisms $f:X\rTo\cX$ with $X\in\cM$.

The morphisms from $f:X\rTo\cX$ to $g:Y\rTo\cX$ are given by 
pairs $(\alpha,\beta)$ where
$\alpha:X\rTo Y$ and $\beta:f\to g\alpha$ is a 2-morphism.
Coverings in $\cX_\et$ are generated by surjective  \'etale 
morphisms. 
We denote by $\cO_{\cX}$ or simply $\cO$ the structure sheaf of $\cX$
of respectively regular / holomorphic / smooth functions.

The sheaves of $\cO$-modules on $\cX$ form a tensor (=symmetric monoidal)
category.

{\bf Notation.} {\em In what follows $\Sh(\cX)$ denotes the category of
sheaves of $\cO$-modules in the complex-analytic or $C^\infty$ case, and
of quasicoherent sheaves in the algebraic case.
}

\subsubsection{Inverse and direct image}
Any
1-morphism of orbifolds $f:\cX\rTo\cY$ induces  a pair of adjoint
functors
$$ f^*:\Sh(\cY) \pile{\rTo \\ \lTo } \Sh(\cX):f_*$$
of inverse and direct image.

In order to define them it is convenient to use a bigger site $\cX_\ET$
whose objects are all \'etale 1-morphisms $f:\cY\rTo\cX$ and
morphisms from $f:\cY\rTo\cX$ to $f':\cY'\rTo\cX$ are given as pairs 
$(\alpha,\beta)$ where
$\alpha:\cY\rTo\cY'$ and $\beta:f\to f'\alpha$ is a 2-morphism.

A 1-morphism of orbifolds $f:\cX\rTo\cY$ induces a morphism
of sites
$$f:\cX_\ET\rTo\cY_\ET$$
which is, by definition, a continuous functor $f^{-1}:\cY_\ET\rTo\cX_\ET$ 
commuting with finite limits. This induces a pair of adjoint functor
of inverse and direct image in a standard way, see~\cite{sga4}, IV.4.

In the case $f$ is representable there is a morphism of sites
$$f:\cX_\et\rTo\cY_\et$$
and one does not need the ``big'' site $\cX_\ET$.

Let now $f,g:\cX\rTo\cY$ be two 1-morphisms and let $\alpha:f\rTo g$
be a 2-morphism. For any $\cU\in\cY_{\ET}$ an isomorphism
$$\alpha:g^{-1}\cU\rTo f^{-1}\cU$$
induces a canonical morphism
$$\alpha:f_*M(\cU)=M(f^{-1}\cU)\rTo M(g^{-1}\cU)=g_*M(\cU).$$
Thus, we have a canonical morphism
\begin{equation}
\label{2-on-dir-im}
\alpha:f_*\rTo g_*,
\end{equation}
as well as
\begin{equation}
\label{2-on-inv-im}
\alpha:g^*\rTo f^*.
\end{equation}

\subsubsection{Projection formula.}
Let $f:\cX\rTo \cY$ be a 1-morphism of orbifolds. One has the following
canonical morphism of bifunctors.

\begin{equation}
\label{proj-for}
f_*(M)\otimes N\rTo f_*(M\otimes f^*(N))
\end{equation}

induced by the adjunction morphism $f^*f_*(M)\rTo M$.

Projection formula states that for certain classes of $f$
 the map ~(\ref{proj-for}) is an isomorphism.

We need a very special case of this.
\begin{prp}Assume $f:\cX\rTo \cY$ is locally a disjoint union of closed
embeddings. Then the canonical morphism~(\ref{proj-for}) is an isomorphism. 
\end{prp}
\begin{proof}The claim reduces to the case of a closed embedding in $\cM$.
In this case $f_*$ identifies $\Sh(\cX)$ with the full subcategory
of $\Sh(\cY)$ consisting of $\cO_{\cY}$-modules annihilated by the ideal
of functions vanishing on $f(\cX)$. In this case the projection formula is 
obvious.
\end{proof}

\section{The functor $\Theta$}
\label{theta}
\subsection{Drinfeld double of a monoidal category}

Recall the construction of a double of a monoidal category
defined in~\cite{js}\footnote{
Joyal and Street call it ``the center of a monoidal category''.
We feel this name is somewhat misleading since $\DD(\cM)\ne\cM$
for a symmetric monoidal category $\cM$.
}.

Recall that a monoidal category $\cM$ is a category endowed
with a bifunctor $M,N\mapsto M\otimes N$ and an isomorphism 
$a:(M\otimes N)\otimes K\rTo M\otimes (N\otimes K)$
(associativity constraint) satisfying the pentagon identity. 
The unit object $\one$ together with isomorphisms 
$\id_{\cM}\rTo \_\otimes\one$,
$\id_{\cM}\rTo \one\otimes\_$, is given, so that the two compositions
from $\id_{\cM}$ to $\one\otimes\_\otimes\one$ coincide.

The associativity constraint determines an isomorphism between any
pair of bracketings of a sequence of objects in $\cM$. These isomorphisms
being compatible, this allows to write $M_1\otimes\ldots\otimes M_n$ for
the tensor product of $n$ objects, omitting the brackets.

Let $\cM$ be a monoidal category. Its double $\DD(\cM)$ is defined as follows.
An object of $\DD(\cM)$ is a pair $(A,\theta)$ where $A\in\cM$ and 
$\theta:A\otimes\_\rTo \_\otimes A$ is an isomorphism of functors
satisfying the two {\em factorization properties}:
\begin{itemize}
\label{factorization}
\item[(F1)] $\theta(\one)=\id_A$.
\item[(F2)] $\theta(B\otimes C)=(B\otimes\theta(C))\circ(\theta(B)\otimes C).$
\end{itemize}
A morphism $(A,\theta)\rTo (A',\theta')$ is a morphism $f:A\rTo A'$
compatible with $\theta$ and $\theta'$. 

The monoidal structure in $\DD(\cM)$ is given by the formula
$$ (A,\theta)\otimes(A',\theta')=(A\otimes A',\theta'')$$
where $\theta'':A\otimes A'\otimes X\rTo X\otimes A\otimes A'$
is defined as the composition of $\theta'$ with $\theta$.

The monoidal category $\DD(\cM)$ is braided: the isomorphism 
$$ (A,\theta)\otimes(A',\theta')\rTo(A',\theta')\otimes(A,\theta)$$
is defined by the map $\theta(A'):A\otimes A'\rTo A'\otimes A$.
\begin{rems}
1. If $\cM$ is a discrete monoidal
category (that is a monoid), the double $\DD\cM$ is the center of the monoid.

2. The double construction for monoidal categories corresponds to
Drinfeld's construction ~\cite{dqg} of the double of a Lie bialgebra.
This is why we use the term {\em Drinfeld double}. 
The construction for monoidal categories is described by 
A.~Joyal and R.~Street in~\cite{js}.
\end{rems}

\subsection{Case of symmetric monoidal category}
In this note we apply the Drinfeld double construction to the tensor
(=symmetric monoidal) category of sheaves on an orbifold.

Let $\cM$ be a tensor category so that, apart from the monoidal structure,
a commutativity constraint $\sigma_{AB}:A\otimes B\rTo^\isom B\otimes A$ is 
given.

An object $(A,\theta)\in\DD(\cM)$ defines uniquely an automorphism
$\tau$ of the functor $M\mapsto A\otimes M$ as the composition of 
$\theta(M)$ with the commutativity constraint. We rewrite below the 
factorization axioms (F1), (F2) in terms of automorphism $\tau$.

\begin{itemize}
\label{factorization-2}
\item[(FT1)] $\tau(\one)=\id_A$.
\item[(FT2)] $\tau(B\otimes C)=
(A\otimes\sigma_{CB})\circ(\tau(C)\otimes B)\circ
(A\otimes\sigma_{BC})\circ(\tau(B)\otimes C)$.
\end{itemize}

\subsection{The canonical functor}
\label{Theta}
Let $\cX$ be an orbifold and let $I\cX$ be its inertia orbifold. 
Now we are ready to define the canonical functor
$$\Theta:\Sh(I\cX)\rTo \DD(\Sh(\cX)).$$
Let $\pi:I\cX\rTo\cX$ be the natural projection.
Let $M\in\Sh(I\cX)$. The object $\Theta(M)\in\DD(\Sh(\cX))$ is given
by the pair $(\pi_*(M),\tau)$ where $\tau$ is the endomorphism of the functor
$N\mapsto \pi_*(M)\otimes N$ defined as follows.

Notice that the projection $\pi:I\cX\rTo\cX$ is endowed with a canonical
2-automorphism $\zeta:\pi\rTo\pi$ assigning to a pair $(x,\gamma)$
the automorphism $\gamma$ of $x\in\cX(M)$.

Thus, the inverse image $\pi^*(N)$ admits a canonical automorphism
$\zeta$ whose fiber at a point $(x,\gamma)\in I\cX_0$ is given by
\begin{equation}
  \label{eq:inv}
  \pi^*(N)_{(x,\gamma)}=N_x\rTo^\gamma N_x=\pi^*(N)_{(x,\gamma)}.  
\end{equation}
By the projection formula
\begin{equation}
  \label{eq:proj-can}
  \pi_*(M)\otimes N=\pi_*(M\otimes \pi^*(N)). 
\end{equation}

We define $\tau:\pi_*(M)\otimes N\rTo\pi_*(M)\otimes N$
as the map induced by the canonical automorphism $\zeta$ of $\pi^*(N)$
described above.

\section{Affine algebraic version}
\label{affine}
It is a pleasure to start with the affine algebraic version of the theorem, 
even though we are unable to deduce from it a global version for DM stacks.

Let $A$ be a commutative ring acted upon by a finite group $G$. Let $\cC$
be the category of $G$-equivariant $A$-modules. We will explicitly
describe now the double $\DD(\cC)$.

\subsection{Calculation}
Let $R=A\rtimes G$ be the twisted group ring, so that $\cC=_R\!\Mod$, the
category of left $R$-modules. This is a monoidal category with respect
to the operation
$$ M,N\mapsto M\otimes_AN$$
defined by the $A$-coalgebra structure on $R$ 
(with respect to the left $A$-module structure on $R$) defined by the formula
 $\Delta(g)=g\otimes g$.

\begin{lem}
\label{Mt}
The functor $N\mapsto M\otimes N$ can be canonically presented
as tensoring by an $R$-bimodule $\widetilde{M}$.
\end{lem}
\begin{proof}
We construct  $\widetilde{M}$ as follows. As a left $R$-module,
$\widetilde{M}$ is just the tensor product $M\otimes R$ with $R$-module
action defined using the coproduct in $R$. The right $R$-module
structure is defined by right multiplication. Thus, one has
\begin{eqnarray}
\label{leftmodule}
ag(m\otimes h)&=&ag(m)\otimes gh\\
\label{rightmodule}
(m\otimes h)ag&=&h(a)m\otimes hg.
\end{eqnarray} 
Notice that $\widetilde{M}$ has different left and right $A$-module structures.
The check that $\widetilde{M}$ is the bimodule we need, is straightforward.
\end{proof}

Note that $\cC$ is a symmetric monoidal category, thus an object of $\DD(\cC)$
is defined by a pair $(M,\tau)$ where $M\in\cC$ and $\tau$ is an endomorphism
of the functor $N\mapsto M\otimes N$ satisfying the 
properties~\ref{factorization-2}.

\begin{lem}
\label{endobim}
Any endomorphism of the functor $N\mapsto M\otimes N$ is uniquely defined
by an $R$-bimodule endomorphism of $\widetilde{M}$.
\end{lem}
\begin{proof}
\end{proof}

Let us describe $R$-bimodule endomorphisms of $\widetilde{M}$.
If $\Phi:M\otimes R\rTo M\otimes R$ is such an endomorphism, one has
$\Phi(m\otimes g)=\phi(m)g$ where $\phi$ is the restriction of $\Phi$
to $M=M\otimes 1\subseteq M\otimes R$.

Let $\phi(m)=\sum_{h\in G}\phi_h(m)\otimes h$.

Then $\Phi(m\otimes g)=\sum_{h\in G}\phi_h(m)\otimes hg$.
The $R$-bilinearity conditions on $\Phi$ amount to the conditions
\begin{eqnarray}
\label{eq:lgphi}
g\phi_h(m)&=&\phi_{ghg^{-1}}(gm)\\
\label{eq:laphi}
\phi_h(am)&=&a\phi_h(m)\\
\label{eq:raphi}
\phi_h(am)&=&h(a)\phi_h(m)
\end{eqnarray}
for all $a\in A,\ g,h\in G,\ m\in M$,
where (\ref{eq:lgphi}) follows from the preservation of left action by 
$g\in G$ and (\ref{eq:laphi}) and (\ref{eq:raphi}) say about preservation 
of the left an the right $A$-module structure respectively.
Comparing the conditions  (\ref{eq:laphi}) and (\ref{eq:raphi}) we deduce
that $\phi_h$ are $A$-module endomorphisms and $h(a)-a$ annihilates $\phi_h$
for each $a\in A$.

Finally, factorization properties (F1), (F2), see~\ref{factorization}, give
the conditions
\begin{eqnarray}
  \label{eq:fact}
  \sum_h\phi_h&=&1\\
  \phi_g\phi_h&=&0,\ g\ne h\\
  \phi_g^2&=&\phi_g.
\end{eqnarray}
Thus the endomorphisms $\phi_g$ are orthogonal idempotents and therefore 
they provide a decomposition of $M$ into a direct sum.
This means that an object of $\DD(\cC)$ is given by an equivariant 
$A$-module 
$$M=\bigoplus_{h\in G}M_h$$
such that $g(M_h)=M_{ghg^{-1}}$ and each $M_h$ is annihilated by $h(a)-a$
for each $a\in A$.

We have proven the following

\begin{thm}
\label{affalg}
Let a finite group $G$ act on a commutative ring $A$. Define $B$
as the direct product
$$ B=\prod_{h\in G}A/\langle h(a)-a\rangle,$$
with the action of $G$ defined by the formula
$$ x=\{x_h\}\in B\Longrightarrow g(x)=\{g(x_{g^{-1}hg})\}.$$
Then the double of the category of $G$-equivariant $A$-module is equivalent
to the category of $G$-equivariant $B$-modules.
\end{thm}

\subsection{}
\label{quasiinv}
A direct calculation shows that the assignment of Theorem~\ref{affalg}
is quasiinverse to the functor $\Theta:\Mod(B)\rTo\DD(\Mod(A))$ defined 
in~\ref{Theta}. 

\section{Affine $C^\infty$ and complex-analytic version}
\label{affine2}
In this section we will check (in a way, similar to that of the previous 
section) that the canonical functor $\Theta$ defined in~\ref{Theta}, 
is an equivalence for affine $C^\infty$ and affine complex-analytic orbifolds.

\subsection{$C^\infty$ case}

Let $X$ be a smooth manifold.
Let $M\in\Sh([G\bs X])$. This is a $G$-equivariant sheaf on $X$.
Define $\widetilde{M}=M\otimes_{\C}\C G=\oplus_{h\in G}M\otimes h$.
This is a $G$-equivariant sheaf on $X$ with the $G$-action given by the
formula
$$ g(m\otimes h)=g(m)\otimes gh,\quad g\in G, m \text{ is a local section of }
M.$$
The sheaf $\widetilde{M}$ admits also a right $G$-action via
$$ (m\otimes h)g=m\otimes hg.$$
Also the global $C^\infty$ functions on $X$ act on $\widetilde{M}$
on the right by the formula
$$ (m\otimes h)a=h(a)m\otimes h.$$
The formulas above are similar to the formulas~(\ref{leftmodule}),
(\ref{rightmodule}).

Any endomorphism of the functor $N\mapsto M\otimes N$ induces
an equivariant endomorphism of the sheaf $\widetilde{M}$ commuting with the
right actions. In fact, let $\cO$ be the sheaf of $C^\infty$ functions on 
$\cX$. Then
$\widetilde{M}=M\otimes\widetilde{\cO}$. Thus, any endomorphism 
$\tau:M\otimes\_\rTo M\otimes\_$ induces an endomorphism  $\Phi=
\tau(\widetilde{\cO})$ of $\widetilde{M}$; it commutes with the right 
actions since right multiplication by $g\in G$ or by $a\in C^\infty(X)$ 
is a morphism of equivariant $\cO$-modules.

The converse is also true, but we are not using this at the moment.
Similarly to the affine algebraic case, we deduce for 
$\Phi(m\otimes g)=\sum\phi_h(m)\otimes hg$ the formulas analogous
to~(\ref{eq:lgphi})--(\ref{eq:raphi}):
\begin{eqnarray}
\label{eq:lgphi2}
g\phi_h(m)&=&\phi_{ghg^{-1}}(gm)\text{ ($m$ is a local section of $M$)}\\
\label{eq:laphi2}
\phi_h(am)&=&a\phi_h(m)\text{ ($a$ and $m$ are local sections of $\cO$ 
and $M$)}\\
\label{eq:raphi2}
\phi_h(am)&=&h(a)\phi_h(m)\text{ ($a$ is a global smooth function) }
\end{eqnarray}
Suppose now that $\tau$ satisfies the factorization properties.
Then $\phi_h$ are orthogonal idempotents so that for $M_g=\phi_g(M)$
$$ M=\bigoplus_{g\in G}M_g.$$
Comparing ~(\ref{eq:laphi2}) to~(\ref{eq:raphi2}) we get that
the sheaf of $\cO_X$-modules $M_g$ is annihilated by any global
$C^\infty$ function of the form $g(a)-a$ where $a$ is a global smooth
function. We claim this implies that $M_g$ is the direct image of
a sheaf on $X^g$ with respect to the closed embedding $X^g\rTo X$.
To prove this we have to check that for each open $U\subseteq X$,
$m\in M(U)$ and $f\in\cO(U)$ such that $f$ is zero on $X^g\cap U$,
one has $fm=0$. This can be checked for arbitrarily small neighborhoods
of points of $X^g$. Since $X^g$ is regularly embedded in $X$, the fact
follows from~\cite{mal}, Thm. VI.1.1.

We have therefore constructed a functor
\begin{equation}
\label{}
\DD(\Sh([G\bs X]))\rTo \Sh(I[G\bs X]).
\end{equation}

Similarly to~\ref{quasiinv}, the functor constructed above is quasi-inverse 
to the one described in~\ref{Theta}.

\subsection{Complex-analytic case}The proof presented above for the 
smooth case works as well in the complex-analytic case if $X$ is a 
Stein manifold; the only nontrivial step of presenting $M_h$ as a direct
image of a sheaf on $X^h$ is due to the fact that the sheaf $\cO$ of 
holomorphic functions on $X$ is generated by its global sections.

\section{Affine orbifold atlas}
\label{sec:atlas}

In this section we describe how to present a complex or $C^\infty$ orbifold 
as a result of gluing of affine orbifolds. 

Notice that, historically, the first definition of orbifold
~\cite{satake} has been given in language of orbifold charts;
in the effective case the passage from Satake definition to the
stack version is described in~\cite{hv}, 2.8, using an appropriate
2-direct limit construction. Here we go in the opposite direction,
starting with an $\cM$-orbifold and presenting it as a direct limit
of affine orbifolds.
This allows to present a sheaf on an orbifold $X$ as a compatible
collection of sheaves on the affine orbifold charts. The construction
of inertia orbifold can be also performed chart-by-chart.

\subsection{Definitions}

\subsubsection{Orbifold charts.}
Let $\cM$ be of one of the types mentioned in the introduction.
An (affine) orbifold chart is a pair $(U,G)$ where $U\in\cM$ and $G$ is a 
finite group acting on $X$ and
\begin{itemize}
\item $X$ is affine in the algebraic case.
\item $X$ is Stein in the complex case.
\item (no extra conditions in the $C^\infty$ case.)
\end{itemize}

A morphism of orbifold charts $(U,G)\rTo(U',G')$ is defined as a pair
of morphisms $\alpha:U\rTo U',\ \gamma:G\rTo G'$ such that
\begin{itemize}
\item they are compatible, i.e. $\alpha(gu)=\gamma(g)\alpha(u)$
for $g\in G,\ u\in U$;
\item The induced map of orbifolds 
$$[\alpha,\gamma]:\left[G\bs U\right]\rTo\left[G'\bs U'\right]$$
is an open embedding.
\end{itemize}

\subsubsection{Atlas.} Let $X$ be an $\cM$-orbifold. Its full atlas
$\Atlas(X)$ is the category defined as follows.
\begin{itemize}
\item The objects are triples $(U,G,\phi)$ where $(U,G)$ is an affine
orbifold atlas and $\phi:[G\bs U]\rTo X$ is an open embedding.
\item A morphism from  $(U,G,\phi)$ to  $(U',G,\phi')$
is given by a triple $(\alpha,\gamma,\theta)$ where 
$(\alpha,\gamma)$ is a morphism of the orbifold charts, and $\theta$
is a 2-morphism from $\phi$ to $\phi'\circ\left[\alpha,\gamma\right]$.
\end{itemize}

A full subcategory $A\subseteq\Atlas(X)$ is called an affine orbifold atlas
of $X$ if the following properties are satisfied.
\begin{itemize}
\item The map $\coprod_{(U,G,\phi)\in A}[G\bs U]\rTo X$ is surjective.
\item For each point $x\in |X|$ lying in the image of two charts
$a_i=(U_i,G_i,\phi_i)\in A,\quad i=1,2$ there exists a chart 
$a_{12}=(U_{12},G_{12},\phi_{12})\in A$ and a pair of maps
$a_{12}\rTo a_i$ such that $x$ belongs to the image of $a_{12}$.
\end{itemize}

Proposition~\ref{prp:orbi-neighborhood} shows that complex and $C^\infty$ 
orbifolds admit affine orbifold atlases. For DM algebraic stacks the
existence of an affine orbifold atlas seems to be a very strong condition
even over a separably closed field.

\subsection{Sheaves}
Let $D$ be a category. Recall that a $D$-category is, by definition, a functor
$\pi:X\rTo D$.

The sheaves on the orbifold charts of an atlas $A$ of $X$ can be organized 
into a category 
$\Sh$ bifibered over $A^\op$. Put $D=A^\op$.

We define the $D$-cateory $\Sh$ as follows.

The objects of $\Sh$ over $(U,G,\phi)\in D$ are the sheaves on $[G\bs U]$.
An arrow $v:F\rTo F'$ over an arrow $u:(U,G,\phi)\rTo(U',G',\phi')$
is just a map of sheaves $u:F\rTo F'$ over a morphism of orbifolds 
$u^\op:[G'\bs U']\rTo\left[G\bs U\right]$. 
Since $u^\op$ is representable, this is nothing but a compatible collection
of maps $F(V)\rTo F'(V\times_{\left[G\bs U\right]}\left[G'\bs U'\right])$.

It is sometimes convenient to add an initial object to $D=\Atlas(X)^\op$
and to complete $\Sh$ with the sheaves on $X$. Here are the appropriate
definitions.

The category $D_+$ is obtained from $D$ by joining an initial object 
$\varnothing$. We define the $D_+$-category $\Sh_+$ as follows.
Its restriction to $D$ is the category $\Sh$ constructed above. The
objects over $\varnothing$ are just the sheaves on $X$. An arrow
from $F$ to $F'$ over the unique morphism $(d):\varnothing\rTo d$
from $\varnothing$ to $d=(U,G,\phi)$ is a map of sheaves $F\to F'$ over the
map $\phi:[G\bs U]\rTo X$.

\subsubsection{Cofibered, fibered and bifibered categories.} 
Let $\pi:X\rTo D$ be a $D$-category.
Recall that a morphism $v:x\rTo x'$ in $X$ over $u:d\rTo d'$ in $D$
is called cocartesian if the natural map
$$ \Hom_{X_{d'}}(x',y)\rTo\Hom_u(x,y)$$
is a bijection. A $D$-category $X$ is called cofibered if the following
properties are satisfied:
\begin{itemize}
\item For each $u:d\rTo d'$ in $D$ and for each $x\in X_d$ there exists
a cocartesian arrow $v:x\rTo x'$ over $u$.
\item Composition of cocartesian arrows is cocartesian.
\end{itemize}

There is a dual notion of cartesian arrow and of fibered category,
see~\cite{sga1}, Expos\'e VI. A $D$-category is bifibered if it is
cofibered and fibered simultaneously.

The $D_+$-category $\Sh+$ constructed above is obviously bifibered;
this is equivalent to the existence of direct and inverse image functors
for sheaves on orbifolds. We will mention an extra property of $\Sh_+$.

\begin{dfn}A $D$-category $X$ will be called {\em special bifibered category}
if it is bifibered and all cartesian arrows in $X$ are necessarily
cocartesian.
\end{dfn}

\begin{lem}$\Sh_+$ is a special bifibered $D_+$-category. 
\end{lem}
\begin{proof}
For any open embedding
$$u:X\rTo Y$$
of orbifolds the composition $u^*\circ u_*$ is isomorphic to identity.
This immediately implies the lemma.
\end{proof}

\subsection{Effective descent}

Let $D$ be a category and let $D_+$ be the category obtained by adding
the initial object $\varnothing$ to $D$. If $X$ is a $D_+$-category,
we denote $X_D=D\times_{D_+}X$. This is a $D$-category.

If $\pi:X\rTo D$ is cofibered, one defines 
$\Gamma^\coc(\pi)$ (or $\Gamma^\coc(X)$) to be the category of cocartesian 
sections of $\pi$.

Let $X$ be cofibered over $D_+$. We say that $X$ admits an effective
descent if the obvious functor
\begin{equation}
\Gamma^\coc(X)\rTo\Gamma^\coc(X_D)
\end{equation}
is an equivalence. Note that the functor $\Gamma^\coc(X)\rTo X_\varnothing$
is equivalence since $\varnothing$ is an initial object of $D_+$.

\begin{lem}
The $D_+$-category $\Sh_+$ has an effective descent.
\end{lem}
\begin{proof}
This is fairly standard: a sheaf on $X$ is uniquely given by a compatible 
collection of sheaves on its affine orbifold charts.
\end{proof}

\section{The main result}
\label{localtoglobal}

\subsection{Relative monoidal categories}

In this subsection we will show how to generalize
to $D$-categories the notions of monoidal category and the construction of 
double.

\subsubsection{Cofibered monoidal categories.}
A monoidal structure on a $D$-category $X$ is given by a functor
$$ \otimes:X\times_DX\rTo X$$
of $D$-categories, together with a unit $\one:D\rTo D$
and the usual constraints and compatibilities.

This means that tensor product is defined for pairs of objects having
the same image in $D$ and the morphism $v\otimes w:x\otimes y\rTo x'\otimes y'$
is defined for a pair of morphisms $v:x\to x',\ w:y\to y'$ having the
same image in $D$.

A monoidal cofibered category $\pi:X\rTo D$ is a monoidal $D$-category
such that
\begin{itemize}
\item It is cofibered;
\item If $v:x\to x'$ and $w:y\to y'$ are cocartesian over the same arrow
$u$ of $D$, then $v\otimes w$ is also cocartesian.
\end{itemize}

\subsubsection{Drinfeld double of a monoidal $D$-category.}
Let $\pi:X\rTo D$ be a monoidal $D$-category. Its Drinfeld double
is a $D$-category $\DD(X)\rTo D$ defined as follows.

The objects of $\DD(X)$ are defined by the formula
$$ \Ob\DD(X)=\coprod_{d\in D}\Ob\DD(X_d).$$
For a pair $(x,\alpha)$ and $(y,\beta)$ of objects of $\DD(X)$
and for a morphism $u:\pi(x)\rTo\pi(y)$ in $D$, we define
$\Hom_u((x,\alpha),(y,\beta))$ as the set of all $f\in\Hom_u(x,y)$ such that
for each $v:z\rTo z'$ over $u$ the diagram
\begin{equation}
\label{rel-double}
\begin{diagram}
x\otimes z & \rTo^\alpha & z\otimes x \\
\dTo^{f\otimes v} & & \dTo^{v\otimes f} \\
y\otimes z'& \rTo^\beta & z'\otimes y
\end{diagram}
\end{equation}
is commutative.

In the case $X$ is cofibered over $D$, it is enough to verify 
the condition~(\ref{rel-double}) for cocartesian $v$ only. 

\subsubsection{The functor of sections.}

If $\pi:X\rTo D$ is cofibered monoidal, $\Gamma^\coc(\pi)$ is monoidal:
for two cocartesian sections $F,G:D\rTo X$ their tensor product is
defined as
$$(F\otimes G)(d)=F(d)\otimes G(d).$$

Below we show that in the special bifibered case the functors 
$\Gamma^\coc$ and $\DD$ commute.

A $D$-category $X$ is called special bifibered monoidal category if
it is special bifibered and monoidal cofibered. 

\begin{lem}
Assume $\pi:X\rTo D$ is a special bifibered monoidal category.
Then the Drinfeld double $\pi:\DD(X)\rTo D$ is cofibered.
\end{lem}
\begin{proof}
Let $u:d\to d'$ be an arrow in $D$, $(x,\alpha)\in\DD(X_d)$. We have
to construct a cocartesian arrow $(x,\alpha)\rTo (x',\alpha')$ in $\DD(X)$.
Choose a cocartesian map $v:x\to x'$ in $X$ over $u$. We claim that
the isomorphism $\alpha':x'\otimes\_\rTo\_\otimes x'$ is defined canonically
by the choice of $v$. In fact, for $y\in X_{d'}$ choose
a cartesian arrow $w:y''\to y$ over $u$. By the assumption of the proposition,
$w$ is cocartesian, so $v\otimes w: x\otimes y''\rTo x'\otimes y$ is
also cocartesian. This implies that the composition
$$ x\otimes y''\rTo^\alpha y''\otimes x\rTo^{w\otimes v} y\otimes x'$$
factors uniquely through $v\otimes w:x\otimes y''\rTo x'\otimes y$. 
This factorization defines an isomorphism of functors
$\alpha':x'\otimes\_\rTo\_\otimes x'$. One checks directly that the
morphism $(x,\alpha)\rTo(x',\alpha')$ is cocartesian.
Thus, $\DD(X)\rTo D$ is cofibered. 
\end{proof}

\begin{prp}
Assume that $\pi:X\rTo D_+$ is a special bifibered monoidal category
admitting an effective descent. Then the Drinfeld double \\
$\DD(\pi):\DD(X)\rTo D_+$ also admits an effective descent. In particular,
the categories $\DD(X_\varnothing)$ and $\Gamma^\coc(\DD(X_D))$
are naturally equivalent.
\end{prp}
\begin{proof}
By the assumption the functors
$$ X_\varnothing\lTo\Gamma^\coc(X)\rTo\Gamma^\coc(X_D)$$
are equivalences of monoidal categories. This yields the equivalences
$$ \DD(X_\varnothing)\lTo\DD(\Gamma^\coc(X))\rTo\DD(\Gamma^\coc(X_D)).$$

The functor
$$ \rho:\Gamma^\coc\DD(X)\rTo\Gamma^\coc\DD(X_D)$$
enters into the diagram
$$
\begin{diagram}
& & \DD(X_\varnothing) & & \\
& \ruTo^{\simeq} & &  \luTo^{\simeq} &\\
\Gamma^\coc\DD(X) &  & \rTo^{\simeq} & & \DD\Gamma^\coc(X) \\
\dTo^\rho & & & & \dTo^{\simeq} \\
\Gamma^\coc\DD(X_D) &  & \rTo^\sigma & & \DD\Gamma^\coc(X_D)
\end{diagram}
$$
so that $\rho$ and $\sigma$ become (essentially) quasi-inverse.
\end{proof}

\subsection{Inertia through the orbifold atlas} 

Let $\cX$ be an orbifold and let
$A\subseteq\Atlas(\cX)$ be an affine orbifold atlas of $\cX$.

For $a=(U,G,\phi)\in A$ we define $Ia\in\Atlas(I\cX)$ as $Ia=(V,G,\psi)$
where $V=\coprod_{g\in G}U^g,\quad \psi=I\phi:
 \left[G\bs V\right]=
I\left[G\bs U\right]
\rTo I\cX$.
\begin{prp}
The collection $\{Ia|a\in A\}$ is an atlas of $I\cX$.
\end{prp}
\begin{proof}Straightforward.
\end{proof}

\section{concluding remarks}
\subsection{Algebraic case}
Deligne-Mumford stacks do not admit in general an affine orbifold atlas.
This is why we were unable to prove the main result for DM stacks.

It seems plausible that the notion of inertia groupoid should be
changed in the algebraic case, see~\cite{agv}. This, however, does not remedy
the lack of the affine orbifold charts.

\subsection{Monoidal structures on $\Sh(I\cX)$}
\label{newbradedstructure}
Note that the equivalence $\Sh(I\cX)\to\DD(\Sh(\cX))$ does not preserve
the tensor structure. For instance, the unit object on the right hand side
is presented by the structure sheaf on the {\em nontwisted} component
of $I\cX$. The new monoidal structure on $\Sh(I\cX)$ can be described
as follows. Let $I_2\cX$ be the double inertia orbifold. The objects
of $I_2\cX(M)$ consist of triples $(x,\gamma_1,\gamma_2)$ with
$x\in\cX(M),\ \gamma_1,\gamma_2\in\Aut(x)$ and the morphisms 
$(x,\gamma_1,\gamma_2)\rTo(x',\gamma'_1,\gamma'_2)$ are given by morphisms
$x\to x'$ commuting with the automorphisms $\gamma$'s.
One has three morphisms $p_{1,2},m:I_2\cX\rTo I\cX$ sending the triple
$(x,\gamma_1,\gamma_2)$ to $(x,\gamma_1)$, $(x,\gamma_2)$ and 
$(x,\gamma_1\gamma_2)$ respectively. Then the new monoidal structure
on $\Sh(I\cX)$ assigns to a pair $M,N\in\Sh(I\cX)$ the sheaf
$$ m_*(p_1^*(M)\otimes p_2^*(N)).$$

It is very tempting to express the cup-product of the stringy $K$-theory
\cite{jkk} through this monoidal structure.


\begin{thebibliography}{MMM}
\bibitem[AGV]{agv}D.~Abramovich, T.~Graber, A.~Vistoli,
Algebraic orbifold quantum product, Orbifolds in mathematics and physics
(Madison, WI, 2001), 1--24, Contemp. Math., 310, 2002. 
\bibitem[DQG]{dqg}V.~Drinfeld, Quantum groups, Proceedings ICM, Berkeley, 1986,
798--820, Amer. Math. Soc. Providence, 1987. 
arXiv math.AG/0310116
\bibitem[HV]{hv} V.~Hinich, A.~Vaintrob, Augmented Teichm\"uller
spaces and orbifolds, manuscript in preparation.
\bibitem[JKK]{jkk} T.J.~Jarvis, R.~Kaufmann,T.~Kimura,Stringy K-theory
and the Chern character, preprint math.AG/0502280.
\bibitem[JS]{js} A.~Joyal, R.~Street, Tortile Yang-Baxter operators
in tensor categories, JPAA 71(1991), 43--51.
\bibitem[LMB]{lmb}G.~Laumon, L.~Moret-Bailly, Champs Alg\'ebriques,
Springer, 199?.
\bibitem[M]{mal} B.~Malgrange, Ideals of differentiable functions,
Oxford University Press, 1966.
\bibitem[MP]{moep} I.~Moerdijk, D.~Pronk, Orbifolds, sheaves and groupoids,
K-theory, {\bf 12}(1997), 3--21.
\bibitem[Sa]{satake} I.~Satake, The Gauss-Bonnet theorem for $V$-manifolds,
J. math. Soc. Japan, 9(1957), 464--492.
\bibitem[SGA1]{sga1} Rev\^etements \'etales et groupe fondamental
(SGA 1), Lecture notes in mathematics, 224 (1971).
\bibitem[SGA4]{sga4}Th\'eories des topos et cohomologie \'etale de sch\'emas
(SGA4), tomes 1,2. Lecture notes in mathematics, 269,270 (1972).
Descente cohomologique, SGA 4, Expos\'e 5-bis,


\end{thebibliography}
\end{document}